\documentclass[12pt,twoside]{article}

\usepackage[T1]{fontenc}
\usepackage[latin1]{inputenc}
\usepackage[english]{babel}
\usepackage[babel]{csquotes}

\usepackage{cite}

\usepackage{amssymb}
\usepackage{amsmath}
\usepackage{amsthm}
\usepackage{latexsym}
\usepackage{graphicx}
\usepackage{mathrsfs}
\usepackage{mathrsfs}
\usepackage{cite}

\usepackage[usenames,dvipsnames]{color}
\definecolor{viola}{rgb}{0.3,0,0.7}
\definecolor{ciclamino}{rgb}{0.5,0,0.5}

\def\pier #1{#1}

\def\luca #1{{\color{blue}#1}}
\def\luca #1{#1}

\def\lucabis #1{{\color{red}#1}}
\def\lucabis #1{#1}

\theoremstyle{plain}
\newtheorem{thm}{Theorem}[section]

\theoremstyle{definition}
\newtheorem{rmk}[thm]{Remark}

\def\eps{\varepsilon}

\def\En{\mathbb{N}}

\def\Ar{\mathbb{R}}

\def\Pi{\mathbb{P}}

\def\beq{\begin{equation}}
\def\eeq{\end{equation}}

\def\rarr{\rightarrow}

\def\l|{\left\|}
\def\r|{\right\|}
\def\L2{L^2(\Omega)}
\def\H1{H^1(\Omega)}
\def\gH1{H^1_{0,\Gamma_b}(\Omega)}
\def\nL2#1{\l|#1\r|_{\L2}}
\def\nH1#1{\l|#1\r|_{\H1}}
\def\rarrw{\rightharpoonup}

\def\gL2{L^2(\Gamma)}
\def\ngL2#1{\l|#1\r|_{\gL2}}
\def\n{\mathcal{N}}
\def\normaX #1{{\mathopen \| #1\mathclose \|}_X}

\def\diver{\mathop{\rm div}\nolimits}

\begin{document}
\begin{center}
{\huge From the viscous Cahn-Hilliard equation\\ 
to a regularized forward-backward\\[0.2cm] 
parabolic equation\/\footnote{{\bf 
Acknowledgments.}\quad\rm The first author gratefully acknowledges some 
financial support from the GNAMPA (Gruppo Nazionale per l'Analisi Matematica, 
la Probabilit\`a e le loro Applicazioni) of INdAM (Istituto 
Nazionale di Alta Matematica) and the IMATI -- C.N.R. Pavia.}}
\\[1cm]
{\large\sc Pierluigi Colli$^{(1)}$}\\
{\normalsize e-mail: {\tt pierluigi.colli@unipv.it}}\\[.25cm]
{\large\sc Luca Scarpa$^{(2)}$}\\
{\normalsize e-mail: {\tt \pier{luca.scarpa.15@ucl.ac.uk}}}\\[.25cm]
$^{(1)}$
{\small Dipartimento di Matematica ``F. Casorati'', Universit\`a di Pavia}\\
{\small via Ferrata 1, 27100 Pavia, Italy}\\[.2cm]
$^{(2)}$
{\small Department of Mathematics, \luca{University College London}}\\
{\small \luca{Gower Street, London WC1E 6BT,} United Kingdom}\\[1.5cm]
\end{center}
{\small
{\bf Abstract:} 
A rigorous proof is given for the convergence of the solutions of a viscous Cahn--Hilliard 
system to the solution of the regularized version of the forward-backward parabolic equation, as the 
coefficient of the diffusive term goes to $0$. Non-homogenous Neumann boundary condition are handled for 
the chemical potential and the subdifferential of a possible non-smooth double-well functional is considered in 
the equation.  An error estimate for the difference of solutions is also proved in a suitable norm and with a
specified rate of convergence. 
\\[.3cm]
\\[.3cm]
}

\pagestyle{myheadings}
\newcommand\testopari{\sc Pierluigi Colli and Luca Scarpa}
\newcommand\testodispari{\sc Asymptotic analysis for a \pier{Cahn--Hilliard} equation}
\markboth{\testodispari}{\testopari}


\thispagestyle{empty}

\section{Introduction}
\setcounter{equation}{0}
\label{intro}
In this paper we perform an asymptotic analysis on the viscous Cahn--Hilliard initial-boundary value problem 
 \begin{gather}
    \label{uno}
   \partial_t y  - \Delta w = 0 \quad \hbox{in } Q:=\Omega \times (0,T) , \\
   \label{due}
    w= \tau \partial_t y - \delta \Delta y + \beta (y) + \lambda (x,t) \pi (y) - g(x,t) \quad \hbox{in }  Q ,\\
    \label{tre}
    \partial_{\bf n} w= h(x,t) , \quad  \partial_{\bf n} y= 0  \quad \hbox{in } \Sigma :=\Gamma
    \times (0,T) ,\\
    \label{quattro} 
    y(\, \cdot \,,  0)=y_{0} \quad \hbox{in } \Omega ,
\end{gather}
as $\delta \searrow 0$, to obtain the regularized diffusion problem in which \eqref{due} and \eqref{tre}
are replaced by the respective equation and condition
 \begin{gather}
   \label{duel} 
    w= \tau \partial_t y + \beta (y) + \lambda (x,t) \pi (y) - g(x,t) \quad \hbox{in } Q ,\\
    \label{trel}
    \partial_{\bf n} w= h(x,t)  \quad \hbox{in } \Sigma .
\end{gather}
Here, $\Omega $ is a smooth bounded domain in $\Ar^3$ with boundary $\Gamma$ and $T>0$ stands for some final time. The variable $y$  denotes some order parameter which may represent the concentration of a phase, and in \eqref{uno}--\eqref{due} the other variable $w$ plays as chemical potential. The data of the system are the viscosity coefficient $\tau >0$, the given functions $\lambda$ and $g$ in $Q$, which depend on the space and time positions as well as the boundary datum $h$ on $\Sigma$,  and the initial value $y_0$ in $\Omega$.  Note that the boundary conditions in \eqref{tre} are both of Neumann type, and the one for the chemical potential is preserved in the limit procedure, as in \eqref{trel}. Quite unusually for Cahn--Hilliard type problems, this condition is more general that the standard homegeneous one, and this makes that the mean value of $y$ is not conserved in time, as the simple integration on \eqref{uno} on $\Omega $ would imply. 

The nonlinearities $\beta $ and $\pi$ have different properties, although both are equal to $0$ in the critical value $0$: $\pi$ is simply Lipschitz continuous, while $\beta$ represents a continuous increasing function that is allowed in our analysis to become a maximal monotone graph in $\Ar \times \Ar$. Of course, in the case of a multivaled graph $\beta$, \eqref{due} and \eqref{duel} should be properly meant as inclusions (instead of equations).   The point is that, in both situations,  $\beta $ is  the subdifferential of a convex and continuous function $\widehat{\beta } : \Ar \to [0,+\infty) $ and, letting $\widehat{\pi } (r) = \int_0^r \pi (s) ds$, $r\in \Ar$, it turns out that the contribution  $\beta (y) + \lambda (x,t) \pi (y)$ results from the \luca{(sub)differentiation} of the following energy functional 
$$ y \in L^2(\Omega) \mapsto  \int_\Omega\luca{\left( \widehat{\beta } (y) + \lambda(\, \cdot\, , t ) \widehat{\pi } (y) \right)} , \quad \hbox{for almost every } \, t \in (0,T). $$
Note that, taking for a moment $\lambda$ constant, the sum  
$$\psi (y) := \widehat{\beta } (y) + \lambda \widehat{\pi } (y)$$ 
may be seen as a generalization of the well-known double-well potential 
$ \psi (y ) = (y^2 - 1)^2/4, $
actually corresponding to the choices  $ \beta (y) = y^3 $, $\lambda = 1$, $\pi (y ) = -  y$. The presence of a coefficient $\lambda$ varying in space and time is certainly of interest in some cases, for instance you can  think to families of control problems in which some linearization has been carried out (cf., e.g.,  \cite{CGS, CGSpure, CGSvisc}). 

Let us now comment on the two problems. If we combine \eqref{uno} and \eqref{duel}, but without the term  $\tau \partial_t y$, we obtain the nonlinear diffusion equation 
\begin{equation}
\label{intr5}
 \partial_t y  - \Delta ( \psi' (y) - g ) = 0 , 
\end{equation}
which can be derived starting from the mass-balance law
\begin{equation}\label{intr6}
\partial_t y  + \diver {\mathbf h}=0,
\end{equation}
where $\mathbf h$ denotes the flux of diffusant and is related to the
the gradient of the chemical potential $w$ according to the Fick law: 
\begin{equation}\label{intr7}
    \mathbf h=-m\nabla w ,\ \quad m\,\mbox{ mobility constant taken 1 here}.
\end{equation}
The equation of state 
\begin{equation}\label{intr8}
  w=\psi'(y) - g
\end{equation}
contains a nonlinear relation between $w$ and $y$, via the derivative $\psi'$ of the
energy function~$\psi$. As $\psi$ is in general the sum of a convex
function and a concave perturbation, $\psi'$ is non-monotone and the resulting equation 
\eqref{intr5} may be ill-posed. 

Then, suitable regularizations are in order, and 
the most celebrated one is the elliptic regularization
\begin{equation}\label{intr9}
 w= - \delta \Delta y + \psi'(y) - g ,
\end{equation}
with $\delta $ positive coefficient related to surface tension; \eqref{intr9} leads to the 
the well-known Cahn--Hilliard system\cite{CahH,EllSh}. 
Other choices have been considered in the literature: in particular, 
Novick-Cohen and Pego\cite{NovicP1991TAMS} dealt with the viscous regularization
\begin{equation}\label{intr10}
 w= \tau \partial_t y + \psi'(y) - g ,
\end{equation}
which has been recently revisited by Tomassetti~\cite{Tomas} (see also the list of references 
in~\cite{Tomas}). 
In fact, the mathematical problem studied in \cite{Tomas} turns out to be a special case of 
\eqref{uno}, \eqref{duel}, \eqref{trel}, \eqref{quattro} and the existence proof in 
\cite{Tomas} is carried out by working with a Faedo--Galerkin scheme.  

The reader can also examine the paper \cite{BCT} in which a 
more general non-smooth regularization of the form
\begin{equation}\label{intr11}
 w= \gamma(\partial_t y) + \psi'(y) - g ,
\end{equation}
$\gamma$ being a maximal monotone and coercive graph, is considered and the related initial--boundary value problem is investigated when taking Dirichlet boundary conditions for the chemical potential.   

On the other hand, a combination of viscous and energetic regularization leads to the so-called viscous Cahn--Hilliard equation (cf.~\eqref{due})
\begin{equation} 
w= \tau \partial_t y - \delta \Delta y + \psi'(y) - g ,
\nonumber
\end{equation} 
derived by Novick-Cohen\cite{Novic1988viscous} and analytically investigated 
in \cite{elliott1996cahn, EllSt}, with numerical aspects treated  in \cite{bai1995viscous}.  
Other recent contributions \cite{CMZ, CF1, CGS, GiMiSchi, GiMiSchi2} 
deal with Cahn--Hilliard or viscous Cahn--Hilliard equations with singular potentials. 
Moreover, and this is important for a comparison with our approach, 
we point out the article \cite{CF2}, in which a wide class of 
evolution equations, but with monotone nonlinearity, is obtained as asymptotic limit 
of special Cahn--Hilliard systems. Also, we mention the papers 
\cite{CFGS1, CGSpure, CGSvisc, HW1, WaNa, ZW} regarding optimal control 
problems for some Cahn--Hilliard systems that possibly include dynamic boundary 
conditions. 

In the present contribution, we first recall a precise well-posedness result 
for the problem \eqref{uno}--\eqref{quattro}. The statements of other results are also contained in Section~\ref{results}. Then, in Section~\ref{convergence} we derive a number of estimates, independent of $\delta$, on the solutions to \eqref{uno}--\eqref{quattro} 
and pass to the limit, as $\delta $ tends to $0$, using properties like compactness, monotonicity and lower semicontinuinity. Thus, we prove the convergence to the solution of \eqref{uno}, \eqref{duel}, \eqref{trel}, \eqref{quattro} and this convergence proof, in our opinion, deserves some interest 
since at some points the result cannot be taken for granted. Of course, we have to make assumptions on the data: concerning the maximal monotone graph $\beta $ in \eqref{due} and \eqref{duel}, we have to require that (see the later~\eqref{data''}) the growth at infinity is
controlled by the one of the related convex function $\widehat{\beta}$. But this looks 
quite reasonable in the framework of \eqref{duel} and less restrictive than the growth of order $p$ assumed in \cite{Tomas}. Actually, by this rigorous asymptotic we give an alterative proof of the existence of solutions to the limit problem \eqref{uno}, \eqref{duel}, \eqref{trel}, \eqref{quattro}. Moreover, for this problem  we show the continuous dependence on the data $y_0, \, g, \, h$, and consequently the uniqueness of the solution, in Section~\ref{cont_dep_lim}. Last but not least, strengthening a bit the assumptions on $\lambda$ and $g$, we are able to deduce an asymptotic error estimate for the difference of the solutions of the two problems: the proof is given in Section~\ref{error_estim}.


\section{Main results}
\setcounter{equation}{0}
\label{results}

In this section, we give some precise \pier{formulation of the problems} and state our results. 
Let us first recall the working framework: 
\begin{gather}
  \label{omega}
  \Omega\subseteq\Ar^3 \, \text{ smooth bounded domain}\,, \quad \Gamma:=\partial\Omega, 
  \quad \pier{T >0 \, \hbox{ final time}},\\
  \label{Qt}
  Q_t:=\Omega\times(0,t)\,, \quad \Sigma_t:=\Gamma\times(0,t) \quad\text{for all } t\in(0,T],\\
  \label{Q}
  Q:=Q_T\,, \quad \Sigma:=\Sigma_T\,.
\end{gather}
\pier{Then,} $\Omega$ is the \pier{spatial} three-dimensional domain, $\Gamma$ \pier{denotes} its 
boundary, while $Q$ and $\Sigma$ represent the \pier{related} spatiotemporal \pier{sets.
In addition, we} \pier{adopt} the notations
\beq
  \label{spaces}
  H:=\L2\,, \quad V:=\H1
\eeq
and we will denote the duality pairing between $V^*$ and $V$ with the symbol $\left<\cdot,\cdot\right>$. As usual, we \pier{make} the identification
$H\cong H^*$, so that $H$ is continuously embedded in $V^*$ in the standard way: for all $u\in H$ and $v\in V$, we have
$\left<u,v\right>=\left(u,v\right)_H$, where $\left(\cdot,\cdot\right)_H$ is the inner product of $H$.

Secondly, let us recall some useful tools which are often used in dealing with the Cahn-Hilliard equations:
the reader can refer to \cite[Sec.~2, pp.~979-980]{CGS}. We introduce
the notation
\beq
  \label{mean}
  \pier{v^*_{\,\Omega}}:=\frac{1}{\left|\Omega\right|}\left<v^*,1\right> 
  \quad\text{for all } \pier{v^*}\in V^*,
\eeq
\pier{which specifies} the mean value of \pier{the elements of} $ V^*$, and we recall that
\beq
  \label{norm_med}
  \exists\, C>0 \text{ : }\l|v\r|_V^2\leq C\left(\l|\nabla v\r|_H^2+\left|v_\Omega\right|^2\right) \quad\text{for all } v\in V\,.
\eeq
Moreover, we define the operator
\beq
  \label{N}
  D(\mathcal{N})=\left\{v^*\in V^*: v^*_\Omega=0\right\}\,, \quad \mathcal{N}:D(\mathcal{N})\rarr\left\{v\in V: v_\Omega=0\right\}
\eeq
by setting $\mathcal{N}v^*$ as the unique solution with null mean value
to the generalized elliptic equation with homogeneous Neumann boundary conditions, \pier{i.e.,}
\beq
  \label{defN}
  \mathcal{N}v^*\in V\,, \quad \left(\mathcal{N}v^*\right)_\Omega=0\,, \quad
  \int_{\Omega}{\nabla\mathcal{N}v^*\cdot\nabla z}=\left<v^*,z\right> \quad\text{for all } z\in V\,.
\eeq
Furthermore, if we set $\l|\,\cdot\,\r|_*:V^*\rarr[0,+\infty)$ as
\beq
  \label{norm*}
  \luca{\l|v^*\r|^2_*:=\l|\nabla\mathcal{N}(v^*-v^*_\Omega)\r|^2_H+\left|v^*_\Omega\right|^2}\,, \quad v^*\in V^*\,,
\eeq
\pier{then} $\l|\,\cdot\,\r|_*$ is a norm on $V^*$, equivalent to the usual norm \luca{$\l|\,\cdot\, \r|_{V^*}$}, which makes $V^*$ a Hilbert space.
Finally, we recall that
\beq
  \label{prop1}
  \left<v^*,\mathcal{N}v^*\right>=\l|v^*\r|_*^2 \quad\text{for all } v^*\in D(\mathcal{N})\,,
\eeq
\beq
  \label{prop1_bis}
  \exists\, C>0 \text{ : } \l|\n v^*\r|_V^2\leq C\l|v^*\r|_*^2 \quad\text{for all } v^*\in D(\n)\,,
\eeq
while for all $v^*\in H^1(0,T;V^*)$ such that $v^*(t)\in D(\mathcal{N})$ for all $t\in[0,T]$ we have that
\beq
  \label{prop2}
  \left<\partial_t v^*(t),\mathcal{N}v^*(t)\right>=\frac{1}{2}\, \pier{\frac{d}{dt} \l|v^*(t) \r|_*^2} \quad\text{for a.e. } t\in(0,T)\,.
\eeq

Now, it is time make some rigorous assumptions on the data. We assume that
\begin{align}
  \label{beta}
  \pier{\widehat{\beta}}:\Ar\rarr \pier{[0,+\infty)} \text{ \pier{is} proper,} \hbox{ convex \pier{and} lower semicontinuous} \nonumber\\ 
  \pier{\hbox{with } \widehat{\beta}(0)=0\,, \hbox{ and } \beta:=\partial\pier{\widehat{\beta}}\, \hbox{ denotes the subdifferential}}
  \end{align}
\begin{equation}
  \label{pi}
  \pi:\Ar\rarr\Ar \text{ \pier{is a Lipschitz continuous function, with }} \pi(0)=0\,,
\end{equation}
so that $\beta:\Ar\rarr2^\Ar$ is a maximal monotone operator with domain $D(\beta)$ \pier{and} satisfying the condition $\beta(0)\ni0$.
\luca{Please notice also that we have $D(\,\widehat{\beta}\,) = \Ar$.}

We are now ready to focus on our problem and present the main results of the paper.
As we have anticipated, the aim of the work is to take the limit as $\delta\searrow0$ 
in problem \eqref{uno}--\eqref{quattro} \pier{to show the convergence of the solutions 
and provide} an \pier{asymptotic} estimate of the error. More specifically, we present now four fundamental results.
The first one ensures that problem \eqref{uno}--\eqref{quattro} is well posed in a certain variational formulation; the second one
is the effective convergence result as $\delta\searrow0$ and provides the rigorous formulation of the limit problem;
the third one is a continuous dependence result for the limit problem and the fourth one contains the \pier{asymptotic} error estimate.

In the following, we assume to work in the setting \eqref{omega}--\luca{\eqref{spaces}} and \eqref{beta}--\eqref{pi}; moreover, \pier{we let
\beq
\hbox{$\tau>0$ be fixed} 
\label{pier0} 
\eeq
and additionally require that
\begin{equation}
\luca{\exists\,C>0 \ : \  \left|\beta^0(r)\right|\leq C\left(1+\left|\pier{\widehat{\beta}}(r)\right|\right) \quad\text{for all } r\in \Ar \,, }  \label{data''}
\end{equation} 
where $\beta^0(r)$ denotes the minimum-norm element of $\beta(r)$, for any $r\in \Ar$.
We note that, for example, all functions with polynomial and first-order exponential growths comply with our assumption \eqref{data''}.
In addition, we point out that \eqref{data''} implies (actually, it is equivalent to) the condition
\beq
  \luca{|s|\leq C\left(1+\left|\pier{\widehat{\beta}}(r)\right|\right)\quad \hbox{for all $r\in\Ar$, $s\in \beta (r)$} }
  \label{pier2}
\eeq
for the same constant $C$, as checked precisely in the next remark.
\begin{rmk}
\label{Crescita}
Such an equivalence property holds 
for a more general growth condition
and in the general setting of Hilbert spaces. Indeed, let $X$ be a Hilbert space,  
$\widehat{\beta} : \pier{X}\to[0,+\infty)$ be convex and~l.s.c.\
(thus continuous since it is everywhere defined), and ${\pier \Psi} :[0,+\infty)\to[0,+\infty)$ be  a continuous function. If $\beta :=\partial\widehat{\beta}$ and, for every $u\in X$, 
$\beta^0 (u)$~is the element of $\beta (u)$ having minimal norm, from 
the condition
\beq
  \bigl\| \beta^0(u) \bigr\|_X \leq {\pier \Psi}  \left(\widehat{\beta}(u)\right)
  \quad \hbox{for every $u\in X$,}
  \label{pier3}
\eeq
it follows that
\beq
  \normaX{{\pier\zeta}} \leq {\pier \Psi}  \left(\widehat{\beta}(u)\right)
  \quad \hbox{for all $u\in X$ and all $\zeta\in \beta (u) $}.
  \label{pier4}
\eeq
Let us check that. If $u\in {\pier X}$, ${\pier\zeta}\in \beta (u)$ and $\eps>0$, 
the monotonicity of $\beta$ implies
$$
  \bigl( \beta^0(u+\eps{\pier\zeta}) - {\pier\zeta} , (u+\eps{\pier\zeta}) - u \bigr) \geq 0 ,
  \quad \hbox{whence} \quad
  \normaX{{\pier\zeta}} \leq \normaX{\beta^0(u+\eps{\pier\zeta})}.
$$
Then, by applying \eqref{pier3} to $u+\eps{\pier\zeta}$, we infer
\beq
  \normaX{{\pier\zeta}} \leq {\pier \Psi}  \left(\widehat{\beta}(u+\eps{\pier\zeta})\right)
  \nonumber
\eeq
and letting $\eps\searrow 0$ we recover \eqref{pier4} thanks to  the continuity of ${\pier \Psi} \circ\widehat{\beta}$.
\end{rmk}
}%

\pier{Now, we recall a well-posedness result for the problem with $\delta>0$. In order to keep a convenient  notation henceforth, differently from \eqref{uno}--\eqref{quattro} now  we put the subscript $\delta$ to the solution and denote by $y_{0,\delta}$ the initial value corresponding to $\delta$}. 

\begin{thm}
  \label{thm1}
  Let $\delta\in(0,1)$ \pier{and assume} that
  \begin{gather}
    \label{lambda_delta}
    \lambda\in L^\infty(Q)\,,\\
    \label{data_delta}
    \pier{g\in L^2(Q)\,, \quad h\in L^2(\Sigma)}\,,\\
    \label{data'_delta}
     \pier{y_{0,\delta}\in\H1 \ \hbox{ and } \  \widehat{\beta}{}}(y_{0,\delta})\in L^1(\Omega)\,.
  \end{gather}
\pier{Then,} there exist
  \begin{gather}
    \label{y_delta}
    y_\delta\in H^1(0,T;\pier{H})\cap L^\infty(0,T;V)
    \cap{\pier{L^2\bigl(0,T;H^2(\Omega)\bigr)}}\\
    \label{w_xi_delta}
    w_\delta\in L^2(0,T;V)\,, \quad \xi_\delta\in L^2(0,T;H)
  \end{gather}
  satisfying for almost all $t\in(0,T)$ the variational \pier{equalities}
  \begin{gather}
    \label{1}
    \left<\partial_ty_\delta(t),v\right>+\int_\Omega{\nabla w_\delta(t)\cdot\nabla v}=\int_{\Gamma}{h(t)v} \qquad\text{for all } v\in\H1\,,
  \end{gather}
    \begin{align}
      \int_\Omega{w_\delta(t)v}=\ &\tau\int_{\Omega}{\partial_ty_\delta(t)v}+\delta\int_\Omega{\nabla y_\delta(t)\cdot\nabla v}
      +\int_\Omega{\xi_\delta(t)v} \nonumber\\
      &+\int_\Omega{\lambda(t)\pi(y_\delta(t))v}-\int_\Omega{g(t)v} \qquad\text{for all } v\in\H1 \label{2}
     \end{align}
  and such that
  \begin{gather}
    \label{incl_delta}
    \xi_\delta\in\beta(y_\delta) \quad\text{a.e.~in } Q\,,\\
    \label{init_delta}
    y_\delta(0)=y_{0,\delta}\,.
  \end{gather}
  Furthermore, if
  \beq
    \label{data_dep_delta}
    \left(y_{0,\delta}^1, g_1, h_1\right), \ \left(y_{0,\delta}^2, g_2, h_2\right)\in
    \H1\times L^2(Q)\times L^2(\Sigma)\,,
  \eeq
  \beq
  \label{data_dep_delta1}
   \pier{ (y_{0,\delta}^1)_\Omega=(y_{0,\delta}^2)_\Omega\,, \quad \int_\Gamma h_1(t)  = \int_\Gamma h_2(t)}  \quad\text{for a.e. } t\in(0,T)
  \eeq
  \pier{and we let}
  \beq
    \left(y_\delta^1, w_\delta^1, \xi_\delta^1\right), \ \left(y_\delta^2, w_\delta^2, \xi_\delta^2\right)
  \eeq
\pier{denote} \luca{any} respective solutions to problem \eqref{1}--\eqref{init_delta} \pier{corresponding to the data in \eqref{data_dep_delta}--\eqref{data_dep_delta1}}, \pier{then}
  there exists a positive constant $C$\pier{, depending only on $\tau$, the constant in \eqref{prop1_bis}, 
$  \| \lambda\|_{L^\infty (Q)},$  a Lipschitz constant for $\pi$, $\Omega$, and $T$,} such that
    \begin{align}
    &\l|y_\delta^1-y_\delta^2\r|_{L^\infty(0,T;V^*)} 
    + \pier{\tau^{1/2}}\l|y_\delta^1-y_\delta^2\r|_{L^\infty(0,T;H)}+
    \delta^{1/2}\l|\nabla(y_\delta^1-y_\delta^2)\r|_{L^2(0,T;H)}\nonumber\\
    &\leq C\left[ \pier{\l|y_{0,\delta}^1-y_{0,\delta}^2\r|_{*} +\tau^{1/2}} \l|y_{0,\delta}^1-y_{0,\delta}^2\r|_{H}+\l|g_1-g_2\r|_{L^2(Q)}+\l|h_1-h_2\r|_{L^2(\Sigma)}\right]\,.
     \label{dep_delta}
    \end{align}
\end{thm}
  
  Please note that the equations \eqref{1}--\eqref{2} are the natural 
  variational formulations of \eqref{uno} and \eqref{due}, obtained 
  testing by $v\in\H1$ and integrating by parts\pier{, on account of the boundary} conditions in \eqref{tre}.

\begin{rmk}
\label{Prova-esist}
The proof of Theorem \ref{thm1} is omitted in this work, since the reader can refer to 
\pier{similar results shown in, e.g., \cite[Theorems~2.2 and~2.3]{CGS} and
\cite{CF2} for related details. For the sake of completeness, the key idea 
is to approximate the problem using the Yosida regularization $\beta_\eps$
instead of $\beta$ and recover a solution of the approximating problem. 
Then, some uniform estimates are found for this family of solutions 
and a passage to the limit as $\eps\searrow0$ provides the solution to 
the original problem. Let us notice that here the approximation in $\eps$ 
should be carried out with $\delta\in(0,1)$ fixed. 
However, we think that the reader can reconstruct the basic steps of 
the proofs by examining the estimates (independent of $\delta$) we will 
prove in Section~\ref{convergence} and compare with Section~\ref{cont_dep_lim}
for the proof of~\eqref{dep_delta}.}
\end{rmk}

\begin{thm}
  \label{thm2}
  \pier{Assume \eqref{lambda_delta}--\eqref{data_delta} and 
  \begin{gather}
    \label{data'}
    \pier{y_0\in\L2 \ \hbox{ and } \  \widehat{\beta}}(y_0)\in L^1(\Omega)\,.
  \end{gather}
Then,} there exist a family
  \beq
    \label{approx_init}
    \{y_{0,\delta}\}_{\delta\in(0,1)} \subseteq\H1\,, \quad y_{0,\delta}\rarr y_0 \text{ in } \L2 \text{ as } \delta\searrow0
  \eeq
  and a positive constant $M$ such that
  \beq
    \label{approx_init'}
    \delta^{1/2}\l|\nabla y_{0,\delta}\r|_H\leq M\,, \quad \l|\pier{\widehat{\beta}}(y_{0,\delta})\r|_{L^1(\Omega)}\leq M
     \qquad\text{for all } \, \pier{\delta\in (0,1)}\,.
  \eeq
  Furthermore, for every $\delta\in(0,1)$ let
  \begin{gather}
    \label{y_delta_conv}
    y_\delta\in H^1(0,T;\pier{H})\cap L^\infty(0,T;V)
    \cap\pier{L^2}\bigl(0,T;H^2(\Omega)\bigr)\\
    \label{w_xi_delta_conv}
    w_\delta\in L^2(0,T;V)\,, \quad \xi_\delta\in L^2(0,T;H)
  \end{gather}
  be the solutions to problem \eqref{1}--\eqref{init_delta} with data $(y_{0,\delta}, g, h)$\pier{.
Then, there exist a triplet
  \begin{gather}
    \label{y}
    \pier{y\in H^1(0,T;H)}\,,\\
    \label{w_xi}
    w\in L^2(0,T;V)\,, \quad \xi\in L^2(0,T;H),
  \end{gather}
which solves the problem
  \begin{gather}
    \label{1_lim}
    \int_{\Omega}{\partial_t y(t)v}+\int_{\Omega}{\nabla w(t)\cdot\nabla v}=\int_\Gamma{h(t)v}
    \quad\text{for all }    v\in\H1 \pier{\text{ and a.e. }} t\in(0,T)\,,
 \end{gather}   
     \begin{align}
      \int_\Omega{w(t)v}=\tau\int_{\Omega}{\partial_ty(t)v}
      +\int_\Omega{\xi(t)v}
      +\int_\Omega{\lambda(t)\pi(y(t))v}-\int_\Omega{g(t)v} \nonumber\\
      \text{for all }v\in\pier{\L2}\pier{\text{ and a.e. }}  t\in(0,T)\,, \label{2_lim}
    \end{align}
 \begin{gather}
     \label{incl}
    \xi\in\beta(y) \quad\text{a.e. in } Q\,,\\
    \label{init}
    y(0)=y_0\,,
\end{gather}
and a subsequence $\{\delta_n\}_{n\in\En}$, with $\delta_n\searrow0$ as $n\rarr\infty$, such that
  \begin{gather}
    \label{conv_y_1}
    y_{\delta_n}\rarrw y \quad\text{in } H^1(0,T;H)\,,\\    
    \label{conv_y_2}
       y_{\delta_n}\rarr y \quad\text{in } C^0([0,T];V^*)\cap 
       L^2(0,T;H)\,,\\
    \label{conv_w}
    w_{\delta_n}\rarrw w \quad\text{in } L^2(0,T;V)\,,\\
    \label{conv_xi}
    \xi_{\delta_n}\rarrw \xi \quad\text{in } L^2(Q). 
  \end{gather}
}%
  \end{thm}

\begin{rmk}
\pier{Theorem~\ref{thm2} is the effective convergence result for our problem. 
Note that the convergence properties \eqref{conv_y_1}--\eqref{conv_xi} hold in principle for a subsequence $\{\delta_n \}$ but the next result we state will entail, in particular, the uniqueness of the solution component $y$. About $w$ and $\xi$, they are not unique in general
but their difference $w-\xi$ is uniquely determined from \eqref{2_lim}. Then, we can at least claim that the convergence of $y_\delta $ to $y$ and of $w_{\delta} - \xi_{\delta}$ to $w-\xi$ is ensured for the entire family as $\delta \searrow 0$.
}
\end{rmk}

\begin{thm}
  \label{thm3}
  \pier{Assume \eqref{lambda_delta} and  
  \begin{gather}
    \label{data_dep_lim}
   \quad g_1, g_2\in L^2(Q)\,, \quad h_1, h_2\in L^2(\Sigma)\,,\\
    \label{data'_dep_lim}
     y_0^1, y_0^2\in\L2\,,  \quad \pier{\widehat{\beta}}(y_0^1),  \, 
      \pier{\widehat{\beta}}(y_0^2) \in L^1(\Omega) \,,\\
    \label{data''_dep_lim}
    (y_0^1)_\Omega=(y_0^2)_\Omega\,, \quad \int_\Gamma h_1(t) = \int_\Gamma h_2(t) \quad\text{for a.e. } t\in(0,T)\,.
  \end{gather}
 Let $(y_1, w_1, \xi_1)$ and $(y_2, w_2, \xi_2)$ be \luca{any}
 respective solutions of the limit problem \eqref{1_lim}--\eqref{init} corresponding to
  the data in \eqref{data_dep_lim}--\eqref{data''_dep_lim}. 
  Then, there exists a positive constant $C$\pier{, depending only on $\tau$, the constant in \eqref{prop1_bis}, 
$  \| \lambda\|_{L^\infty (Q)},$  a Lipschitz constant for $\pi$, $\Omega$, and $T$,}  such that}
  \beq
    \label{dep_lim}
    \l|y_1-y_2\r|_{L^\infty(0,T;H)}\leq C\left[\l|y_0^1-y_0^2\r|_H+\l|g_1-g_2\r|_{L^2(Q)}+\l|h_1-h_2\r|_{L^2(\Sigma)}\right]\,.
  \eeq
\end{thm}

\begin{rmk}
  Please note that hypothesis \eqref{data''_dep_lim} is the natural generalization that takes place when 
  dealing
  with a non-homogeneous Neumann boundary condition. As a matter of fact, in the case of homogeneous 
  Neumann conditions for $y$,
  the natural requirement is that $y_0^1$ and $y_0^2$ have the same mean value (see for example 
  \cite{CF1, CGS});
  \pier{when a boundary datum} is introduced, we need to require also that $h_1$ and $h_2$ 
  have the same mean value on $\Gamma$, in order to
  recover two solutions $y_1$ and $y_2$ with same mean value on $\Omega$ \pier{at almost every time, 
  so} allowing us to prove the continuous dependence result.
\end{rmk}

\begin{thm}
  \label{thm4} \pier{Assume that
  \begin{gather}
    \label{lambda-reg}
    \lambda\in L^\infty(Q)\cap L^2\bigl(0,T;W^{1,\infty}(\Omega)\bigr)\,,\\
    \label{data-reg}
    g\in L^2(0,T;V)\,, \quad h\in L^2(\Sigma)\,,
  \end{gather}
besides \eqref{data'}. Then,}
  there exists a positive constant $C$ such that
  the following \pier{asymptotic} estimate holds for all $\delta\in(0,1)$:
  \beq
    \label{err}
    \l|y-y_\delta\r|_{L^\infty(0,T;H)}\leq C\left[\delta^{1/4}+\l|y_0-y_{0,\delta}\r|_H  
    \right]\,.
  \eeq
\end{thm}


\section{The convergence result}
\setcounter{equation}{0}
\label{convergence}

In this section, we present the proof of the convergence result contained in Theorem \ref{thm2}. In particular, we will firstly check 
that an approximation on the initial data satisfying 
\pier{\eqref{approx_init}--\eqref{approx_init'}} actually exists; then, we will find some 
uniform estimates on the solutions $(y_\delta, w_\delta, \xi_\delta)$ and pass to the limit as $\delta\searrow0$.

Let us \pier{specify} some useful notation that we use in \pier{the sequel}. If we
test equation \eqref{1} by $v=\pier{1/|\Omega|}$\pier{,} we deduce that
\beq
  \left(\partial_t y_\delta(t)\right)_\Omega= \frac{1}{|\Omega|}\int_\Gamma{h(t)} \quad\text{for a.e. } t\in(0,T)\,,
\eeq
and, in view of \eqref{init_delta}, 
\beq
  \left(y_\delta(t)\right)_\Omega=\left(y_{0,\delta}\right)_\Omega+\frac{1}{|\Omega|}\int_0^t\!\!{\int_\Gamma{h(s)}\,ds} \quad\text{for all } t\in[0,T]\,.
\eeq
Hence, it is natural to introduce 
\beq
  \label{M_delta}
  M_\delta(t):=\left(y_{0,\delta}\right)_\Omega+\frac{1}{|\Omega|}\int_0^t\!\!{\int_\Gamma{h(s)}\,ds}\,,
\eeq
so that
\beq
  \label{prop_M_delta}
  M_\delta= \left(y_\delta\right)_\Omega\,, \quad M_\delta'=\left(\partial_t y_\delta\right)_\Omega\,.
\eeq
\pier{Owing to \eqref{data_delta} and \eqref{data'} it turns out that
\beq 
  M_\delta  \hbox{ is bounded in } H^1(0,T) \hbox{ independently of }   \pier{\delta\in (0,1)}.
  \label{pier5}
\eeq}%

\subsection{Existence of an \pier{approximating family} $\{y_{0,\delta}\}$}
\label{exist_approx}

For every $\delta\in(0,1)$, let us define $y_{0,\delta}$ as the solution to the \pier{elliptic problem}
\beq
  \label{ellip}
  \begin{cases}
    y_{0,\delta}-\delta\Delta y_{0,\delta}=y_0 \quad\text{in } \Omega\,,\\[0.2cm]
    \pier{\partial_\nu y_{0,\delta}}=0 \quad\text{on } \Gamma\,.
  \end{cases}
\eeq
\pier{It is well known that $y_{0,\delta}\in\pier{H^2(\Omega)}$ 
and it} satisfies the variational equation
\beq
  \label{ellip'}
  \int_{\Omega}{y_{0,\delta}z}+\delta\int_\Omega{\nabla y_{0,\delta}\cdot\nabla z}=\int_{\Omega}{y_0z} \quad\text{for all } z\in\H1\,;
\eeq
hence, testing \eqref{ellip'} by $z=y_{0,\delta}$, owing to the Young inequality it is easy to see that
\beq
\label{pier1}
  \frac{1}{2}\l|y_{0,\delta}\r|_H^2+\delta\l|\nabla y_{0,\delta}\r|_H^2\leq\frac{1}{2}\l|y_0\r|_{\pier H}^2 \quad\text{for all } \delta>0\,,
\eeq
so that the first estimate in \eqref{approx_init'} is satisfied. Moreover, \pier{from \eqref{pier1} it follows} that there exists $\widetilde{y_0}\in\L2$ 
and a subsequence $\{y_{0,\delta_k}\}_{k\in\En}$ such that
\beq
  \label{y0_tilde}
  y_{0,\delta_k}\rarrw \widetilde{y_0} \quad\text{in } \pier{H}\,, \quad \delta_ky_{0,\delta_k}\rarr0 \quad\text{in } \pier{V} \quad\text{as } k\rarr\infty\,,
\eeq
and letting $k\rarr\infty$ in \eqref{ellip'} we reach
\[
  \int_\Omega{\widetilde{y_0}z}=\int_\Omega{y_0z} \quad\text{for all } z\in \pier{V}\,:
\]
since $V$ is dense in $H$, it turns out that $\widetilde{y_0}=y_0$. Furthermore, the identification of the weak limit implies that
the entire family $\{y_{0,\delta}\}$ weakly converges to $y_0$ in $H$; finally, from \eqref{ellip'} we have
\[
  \limsup_{\delta\searrow0}{\l|y_{0,\delta}\r|_H^2}\leq\l|y_0\r|_H^2\,,
\]
so that $y_{0,\delta}\rarr y_0$ in $H$ and also condition \eqref{approx_init} is satisfied.

It remains to check the second estimate of \eqref{approx_init'}. For every $\eps\in(0,1)$, let $\beta_\eps$ be the Yosida
approximation of $\beta$: hence, since $\beta_\eps$ is Lipschitz continuous,
$\beta_\eps(y_{0,\delta})\in\H1$ and we can test \eqref{ellip'} by $z=\beta_\eps(y_{0,\delta})$, obtaining
\[
  \int_\Omega{\beta_\eps(y_{0,\delta})(y_{0,\delta}-y_0)}+\delta\int_\Omega{\beta'_\eps(y_{0,\delta})\left|\nabla y_{0,\delta}\right|^2}=0\,.
\]
Hence, thanks to the subdifferential \pier{property} and the monotonicity of \pier{$\beta_\eps$} we have
\[
  \int_\Omega{\pier{\widehat{\beta}}_\eps(y_{0,\delta})}-\int_\Omega{\pier{\widehat{\beta}}_\eps(y_0)}\leq
  \int_\Omega{\beta_\eps(y_{0,\delta})(y_{0,\delta}-y_0)}\leq0\,\pier{,}
\]
whence \pier{(cf., e.g., \cite[Thm.~2.2, p.~57]{Barbu})}
\[
  \int_\Omega{\pier{\widehat{\beta}}_\eps(y_{0,\delta})}\leq\int_\Omega{\pier{\widehat{\beta}}_\eps(y_{0})}\leq\int_\Omega{\pier{\widehat{\beta}}(y_0)}
\]
for all $\eps\in(0,1)$. Taking the limit as $\eps\pier{\searrow} 0$ in the \pier{above inequality}, thanks to the Fatou lemma and condition \eqref{data'}
we obtain also the second part of \eqref{approx_init'}.

\subsection{The estimate on $y_\delta$}

We want to find some uniform estimates on $y_\delta$: firstly, let us notice that equation \eqref{prop_M_delta} ensures that\pier{,} 
for almost every $t\in(0,T)$, $y_\delta(t)-M_\delta(t)\in D(\mathcal{N})$, so that $\mathcal{N}(y_\delta(t)-M_\delta(t))$ makes sense.
Hence, we can test equation \eqref{1} with $\mathcal{N}(y_\delta(t)-M_\delta(t))$ and \eqref{2} with $-(y_\delta(t)-M_\delta(t))$:
summing up the two equations, the second and third integral on the left hand side cancel
thanks to \eqref{defN}\pier{. Then,} we obtain
\[
  \begin{aligned}
  &\left<\partial_t(y_\delta(t)-M_\delta(t)),\n(y_\delta(t)-M_\delta(t))\right>+
  \frac{\tau}{2}\frac{d}{dt}\l|y_\delta(t)-M_\delta(t)\r|_H^2+
  \delta\l|\nabla y_\delta(t)\r|^2_H
  \\
  &+\int_\Omega{\xi_\delta(t)\left(y_\delta(t)-M_\delta(t)\right)}=
 -\left<M'_\delta(t),\n(y_\delta(t)-M_\delta(t))\right>+ 
  \int_\Gamma{h(t)\n(y_\delta(t)-M_\delta(t))}
  \\
  &
  -\tau\int_\Omega{M_\delta'(t)(y_\delta(t)-M_\delta(t))}
  -\int_\Omega{\lambda(t)\left(\pi(y_\delta(t))-\pi(M_\delta(t))\right)(y_\delta(t)-M_\delta(t))}
  \\
  &
  -\int_\Omega{\lambda(t)\pi(M_\delta(t))(y_\delta(t)-M_\delta(t))}
  +\int_\Omega{g(t)(y_\delta(t)-M_\delta(t))}\,.
  \end{aligned}
\]
\pier{Recalling \eqref{N}--\eqref{defN}, we note that}
\begin{gather}
  -\left<M'_\delta(t),\n(y_\delta(t)-M_\delta(t))\right>=-M'_\delta(t)|\Omega|\left(\n(y_\delta(t)-M_\delta(t))\right)_\Omega=0\,,
  \nonumber\\
  \pier{{}-\tau\int_\Omega{M_\delta'(t)(y_\delta(t)-M_\delta(t))} = {}
  - \tau M'_\delta(t)|\Omega|\left(y_\delta(t)-M_\delta(t)\right)_\Omega  = 0}
  \nonumber
\end {gather} 
while the subdifferential rule for $\beta$ together with \eqref{incl_delta} leads to
\[
  \int_\Omega{\xi_\delta(t)\left(y_\delta(t)-M_\delta(t)\right)}\geq
  \int_\Omega{\pier{\widehat{\beta}}(y_\delta(t))}-\int_\Omega{\pier{\widehat{\beta}}(M_\delta(t))}\,.
\]
Hence, taking now into account conditions \eqref{prop1_bis}--\eqref{prop2} and using the assumptions \eqref{lambda_delta}--\eqref{data_delta} and \eqref{beta}--\eqref{pi} on the data,
\pier{we integrate on $(0,t)$ and, with the help of} the Young inequality, we deduce that
\[
  \begin{aligned}
    &\frac{1}{2}\l|y_\delta(t)-M_\delta(t)\r|_*^2+\frac{\tau}{2}\l|y_\delta(t)-M_\delta(t)\r|_H^2+\delta\int_0^t\l|\nabla y_\delta(s)\r|_H^2\,ds
    +\int_0^t\!\!\int_\Omega{\pier{\widehat{\beta}}(y_\delta(s))\,ds}\\
    &\leq\frac{1}{2}\l|y_{0,\delta}-(y_{0,\delta})_\Omega\r|_*^2
    +\frac{\tau}{2}\l|y_{0,\delta}-(y_{0,\delta})_\Omega\r|_H^2
    +\frac{1}{2}\l|h\r|_{L^2(\Sigma)}^2\\
    &\pier{{}+\frac{C}{2}\int_0^t\l|y_\delta(s)-M_\delta(s)\r|_*^2\,ds+ \pier{\left(\luca{C_\pi\l|\lambda\r|_{L^\infty(Q)}}+1 \right)} \int_0^t{\l|y_\delta(s)-M_\delta(s)\r|_H^2\,ds}}\\
    &+\frac{1}{2}\pier{\left(C_\pi^2|\Omega|\l|\lambda\r|_{L^\infty(Q)}^2\l|M_\delta
     \r|_{L^2(0,T)}^2 +\l|g\r|^2_{L^2(Q)}\right)}  +\int_0^t\!\!\int_\Omega{\pier{\widehat{\beta}}(M_\delta(s))\,ds}
  \end{aligned}
\]
\pier{for all $t\in(0,T)$ and for some Lipschitz constant $C_\pi$ of $\pi$}.
Now, \pier{in view of  \eqref{pier5} and the continuity of $\widehat{\beta} $ on $\Ar$ (ensured by 
\eqref{beta} and \eqref{data''}) we have
\beq
  \label{prop_M_delta'}
  \l| M_\delta \r|_{L^\infty(0,T)}
  + 
   \l| \widehat{\beta} (M_\delta) \r|_{L^\infty(0,T)}
   + \l| M'_\delta \r|_{L^2(0,T)}
\leq C
\eeq
for some constant $C$ independent of $\delta \in (0,1)$.}
Hence, taking these remarks into account \pier{and 
recalling that $\{y_{0,\delta}\}$ converges in $H\subseteq V^*$ 
by \eqref{approx_init}, we infer that} 
\[
  \begin{split}
    \frac{1}{2}\l|y_\delta(t)-M_\delta(t)\r|_*^2&+\frac{\tau}{2}\l|y_\delta(t)-M_\delta(t)\r|_H^2+\delta\int_0^t\l|\nabla y_\delta(s)\r|_H^2\,ds
    +\int_0^t\!\!\int_\Omega{\pier{\widehat{\beta}}(y_\delta(s))\,ds}\\
    &\leq \pier{C\left(1+ \int_0^t\l|y_\delta(s)-M_\delta(s)\r|_*^2\,ds+
    \int_0^t{\l|y_\delta(s)-M_\delta(s)\r|_H^2\,ds}\right).}
  \end{split}
\]
Then, the Gronwall lemma ensures that \pier{%
\begin{gather}
  \label{est1_provv}
  \l|y_\delta-M_\delta\r|_{L^\infty(0,T;V^*)}+
  \l|y_\delta-M_\delta\r|_{L^\infty(0,T;H)}+
  \delta^{1/2}\l|\nabla y_\delta\r|_{L^2(0,T;H)}\leq C,\\
  \label{est2}
  \l|\pier{\widehat{\beta}}(y_\delta)\r|_{L^1(Q)}\leq C 
\end{gather}
for all $\delta \in (0,1);$
finally, since $y_\delta  = (y_\delta-M_\delta) + M_\delta $, on account of 
\eqref{prop_M_delta'} and \eqref{est1_provv} we find out that}%
\beq
  \label{est1}
  \l|y_\delta\r|_{L^\infty(0,T;V^*)}+\l|y_\delta\r|_{L^\infty(0,T;H)} +
  \delta^{1/2}\l|y_\delta\r|_{L^2(0,T;V)}\leq C \quad\text{for all } \delta\in(0,1)\,.
\eeq

\subsection{The estimate on $\partial_t y_\delta$}

Let us now prove some uniform \pier{estimate} on $\partial_ty_\delta$. 
\pier{Observe that, in view of the regularity \eqref{y_delta}, the variational equality \eqref{2} yields the equation 
\beq
\label{pier6}
 w_\delta= \tau\partial_t y_\delta- \delta\Delta y_\delta + \xi_\delta + 
 \lambda \pi(y_\delta) - g  \quad\text{a.e. in } Q
\eeq
along with the Neumann homogeneous boundary condition $\partial_\nu y_\delta =0$ a.e. on 
$\Sigma$. Hence, we can take $v=\n(\partial_ty_\delta(t)-M'_\delta(t))$ in \eqref{1} and test
\eqref{pier6} at time $t$ by $-(\partial_ty_\delta-M'_\delta)(t)$. Note that this makes sense since  $-(\partial_ty_\delta-M'_\delta)$ is in $ L^2 (Q)$.}
Summing up the two equations, the second and third integral on the left hand side cancel thanks to \eqref{defN}: hence, using condition \eqref{prop1}
and integrating on $(0,t)$,
for almost every $t\in(0,T)$ we \pier{have} 
\[
  \begin{aligned}
  &\int_0^t{\l|\partial_ty_\delta(s)-M'_\delta(s)\r|_*^2\,ds}+\tau\int_0^t{\l|\partial_ty_\delta(s)-M'_\delta(s)\r|_H^2\,ds}
  +\frac{\delta}{2}\l|\nabla y_\delta(t)\r|_H^2
  \\
  &+\int_0^t\!\!\int_\Omega{\xi_\delta(s)\partial_ty_\delta(s)\,ds}=
  \frac{\delta}{2}\l|\nabla y_{0,\delta}\r|_H^2
  \pier{{}-\int_0^t\left<M'_\delta(s),\n(\partial_ty_\delta(s)-M'_\delta(s))\right>\,ds} 
  \\
  &\pier{{}+\int_0^t\!\!\int_\Gamma{h(s)\n(\partial_ty_\delta(s)-M'_\delta(s))\,ds}}
    -\tau\int_0^t\!\!\int_\Omega{M'_\delta(s)(\partial_ty_\delta(s)-M'_\delta(s))\,ds}
  \\
  &+\int_0^t\!\!\int_\Omega{\xi_\delta(s)M'_\delta(s)\,ds}
  +\int_0^t\!\!\int_\Omega{\left[g(s)-\lambda(s)\pi(y_\delta(s))\right](\partial_ty_\delta(s)-M'_\delta(s))\,ds}\,.
  \end{aligned}
\]
As in the previous \pier{subsection}, using \pier{\eqref{N}--\eqref{defN}} we see that
\[
 \pier{{} -\left<M'_\delta(s),\n(y_\delta(s)-M_\delta(s))\right>=0\,, \quad \int_\Omega{M'_\delta(s)(\partial_ty_\delta(s)-M'_\delta(s))\,ds}=0 \, ,}
\]
while a well-known result (contained for example in \pier{\cite[p.~73]{Brezis}}) ensures that
\[
  \int_\Omega{\xi_\delta(t)\partial_ty_\delta(t)}=\frac{d}{dt}\int_\Omega\pier{\widehat{\beta}}(y_\delta(t))
\]
hence, using the Young inequality and taking into account conditions \eqref{prop1_bis}, \eqref{beta}--\eqref{pi}, \eqref{lambda_delta}--\eqref{data_delta}
and the growth assumption \pier{\eqref{pier2} (cf.~\eqref{data''} and Remark~\ref{Crescita})},
we deduce that for almost every $t\in(0,T)$
\[
  \begin{aligned}
    &\int_0^t{\l|\partial_ty_\delta(s)-M'_\delta(s)\r|_*^2\,ds}+\tau\int_0^t{\l|\partial_ty_\delta(s)-M'_\delta(s)\r|_H^2\,ds}
    +\frac{\delta}{2}\l|\nabla y_\delta(t)\r|_H^2+\int_\Omega{\pier{\widehat{\beta}}(y_\delta(t))}\\
    &\leq\frac{\delta}{2}\l|\nabla y_{0,\delta}\r|_H^2+\int_\Omega{\pier{\widehat{\beta}}(y_{0,\delta})}
    \pier{{}+\frac{1}{2}\int_0^t\l|\partial_ty_\delta(s)-M'_\delta(s)\r|_*^2\,ds
    +\frac{C}{2}\l|h\r|_{L^2(\Sigma)}^2}\\
    &\pier{{}+\frac{\tau}{2}\int_0^t{\l|\partial_ty_\delta(s)-M'_\delta(s)\r|_H^2\,ds}
    +\frac1{2\tau} \l|g-\lambda\pi(y_\delta)\r|^2_{L^2(Q)}}\\
    &+C\int_0^t{M'_\delta(s)\int_\Omega\pier{\widehat{\beta}}(y_\delta(s))\,ds}+C|\Omega|\l|M'_\delta\r|_{L^1(0,T)}\,.
  \end{aligned}
\]
for \pier{some constant} $C>0$. Now, thanks to conditions \eqref{approx_init'}, \pier{\eqref{est1}, \eqref{prop_M_delta'}} and
the Lipschitz continuity of $\pi$, we can \pier{apply the Gronwall lemma and infer  that 
(updating the value of $C$, as usual)
\begin{gather}
  \l|\partial_ty_\delta-M'_\delta\r|_{L^2(0,T;V^*)}+\l|\partial_ty_\delta-M'_\delta\r|_{L^2(0,T;H)}
  +\delta^{1/2}\l|\nabla y_\delta\r|_{L^\infty(0,T;H)}  \leq C, \label{est3}\\
  \l|\pier{\widehat{\beta}}(y_\delta)\r|_{L^\infty(0,T;L^1(\Omega))}\leq C
  \label{est3_bis}
\end{gather}
for all $\delta\in(0,1)$. Hence,  using \eqref{prop_M_delta'} and \eqref{est1} again 
we conclude that}
\begin{gather}
  \label{est4}
\pier{  \l|y_\delta\r|_{H^1(0,T;H)} + \delta^{1/2}\l|y_\delta\r|_{L^\infty(0,T;V)}}\leq C \quad\text{for all } \delta\in(0,1)\,.
\end{gather}

 \subsection{The estimate on $w_\delta$}
 
We now \pier{take  $v=w_\delta(t)-(w_\delta(t))_\Omega\in V$ in equation \eqref{1}; by \eqref{norm_med} we obtain}
\[
  \begin{split}
  \int_\Omega\left|\nabla w_\delta(t)\right|^2&=\int_\Omega\left|\nabla(w_\delta(t)-(w_\delta(t))_\Omega)\right|^2\\
  &\leq C\l|h(t)\r|_{L^2(\Gamma)}\l|w_\delta(t)-(w_\delta(t))_\Omega\r|_V
  +\l|\partial_ty_\delta(t)\r|_{V^*}\l|w_\delta(t)-(w_\delta(t))_\Omega\r|_V\\
  &\leq
  C'\left(\l|h(t)\r|_{L^2(\Gamma)}^2+\l|\partial_ty_\delta(t)\r|_{V^*}^2\right)
  +\frac{1}{2}\int_\Omega\left|\nabla w_\delta(t)\right|^2
  \end{split}
\]
\pier{for almost all $t\in(0,T)$ and}
for some two positive constants $C$ and $C'$. Then, recalling the estimate \eqref{est4} just proved and hypothesis \eqref{data_delta} on $h$, we deduce that
\beq
  \label{est6_provv}
  \l|\nabla w_\delta\r|_{L^2(0,T:H)}\leq C \quad\text{for all } \delta\in(0,1)\,.
\eeq
Moreover, \pier{taking  $v=1$ in \eqref{2}, thanks to \eqref{data''} for almost all $t\in(0,T)$ we have
\[
  \begin{split}
 |\Omega|  \left(w_\delta(t)\right)_\Omega&\leq
  \tau\l|\partial_ty_\delta(t)\r|_{L^1(\Omega)}+C\left( |\Omega|  +\l|\pier{\widehat{\beta}}(y_\delta(t))\r|_{L^1(\Omega)}\right)\\
  &+\l|\lambda\r|_{L^\infty(Q)}C_\pi\l|y_\delta(t)\r|_{L^1(\Omega)}+\l|g(t)\r|_{L^1(\Omega)}
  \end{split}
\]
so} that recalling \eqref{est3_bis}, \eqref{est4}, and the hypothesis \eqref{data_delta} on $g$ we \pier{infer}
\beq
\nonumber
  \l|\left(w_\delta\right)_\Omega\r|_{L^2(0,T)}\leq C \quad\text{for all } \delta\in(0,1)\,,
\eeq
\pier{and this uniform bound, along} with \eqref{est6_provv}, ensures that
\beq
  \label{est6}
  \l|w_\delta\r|_{L^2(0,T:V)}\leq C \quad\text{for all } \delta\in(0,1)\,.
\eeq

\subsection{The estimate on $\xi_\delta$}
\label{xi_est}

\pier{Next, we would like to show a uniform estimate for $\xi_\delta$ in $L^2(Q)$. Let us deduce it on some approximating problem in which $\beta$ is replaced by its Yosida regularization $\beta_{\pier{\eps}}
$. Indeed, the estimate, independent of $\eps \in (0,1)$, will be proved for $\beta_{\pier{\eps}}(y_{\delta,\eps})$, where $ y_{\delta,\eps}$ denotes the main component of the approximating solution. Then, 
passing to the limit as $\eps \searrow 0$, one obtains the same estimate for $\xi_\delta$ (cf.~Remark~\ref{Prova-esist}).}

\pier{Let us rewrite \eqref{2} in tems of $w_{\delta,\eps}$ and $ y_{\delta,\eps}$ obtaining}
\[
  \begin{split}
      \int_\Omega{w_{\delta,\eps}(t)v}=&\tau\int_{\Omega}{\partial_ty_{\delta,\eps}(t)v}+\delta\int_\Omega{\nabla y_{\delta,\eps}(t)\cdot\nabla v}
      +\int_\Omega{\beta_\eps(y_{\delta,\eps}(t)v}\\
      &+\int_\Omega{\lambda(t)\pi(y_{\delta,\eps}(t))v}-\int_\Omega{g(t)v} \qquad\text{for all } v\in\H1\,.
   \end{split}
\]
\pier{Taking $v=\beta_{\pier{\eps}}(y_{\pier{\delta, \eps}}(t))\in V$, we have 
\[
  \begin{split}
  &\delta\int_\Omega{\beta'_\delta(y_{\pier{\delta, \eps}}(t))\left|\nabla y_{\pier{\delta, \eps}}(t)\right|^2}+\int_\Omega{\left|\beta_{\pier{\eps}}(y_{\pier{\delta, \eps}}(t))\right|^2}\\
  &=\int_\Omega\left[g(t)-\lambda(t)\pi(y_{\pier{\delta, \eps}}(t))+w_{\pier{\delta, \eps}}(t)-\tau\partial_ty_{\pier{\delta, \eps}}(t)\right]\beta_{\pier{\eps}}(y_{\pier{\delta, \eps}}(t))
  \end{split}
\]
for almost every $t\in (0,T)$. 
Hence, using the monotonicity of $\beta_{\pier{\eps}}$, the hypotheses on $g,\, \lambda$ and $\pi$ together with conditions \eqref{est4} and \eqref{est6}, from integration with respect to time and the elementary Young inequality it follows that
\begin{align*}
  \int_0^T\!\!\int_\Omega\left|\beta_{\pier{\eps}}(y_{\pier{\delta, \eps}})\right|^2
  &\leq   \int_0^T\!\!\int_\Omega \left| g-\lambda\pi(y_{\pier{\delta, \eps}})+w_{\pier{\delta, \eps}}-\tau\partial_ty_{\pier{\delta, \eps}}\right| \, |\beta_{\pier{\eps}}(y_{\pier{\delta, \eps}})| \\
 &\leq C+\frac{1}{2 }\int_0^T\!\!\int_\Omega\left|\beta_{\pier{\eps}}(y_{\pier{\delta, \eps}})\right|^2
\end{align*}
for a constant $C>0$ independent of both $\delta $ and $ \eps \in (0,1)$. Consequently, bearing in mind what we have anticipated before, we conclude that
\beq
  \label{est7}
  \l|\xi_\delta\r|_{L^2(Q)} \leq \liminf_{\eps \searrow 0 } \l| \beta_{\pier{\eps}}(y_{\pier{\delta, \eps}})\r|_{L^2(Q)}
  \leq C \quad\text{for all } \delta\in(0,1)\,.
\eeq
Let us point out a consequence of \eqref{est7}: by a comparison of the terms in \eqref{pier6} 
we deduce that
\beq
  \label{est8_provv}
  \delta\l|\Delta y_\delta\r|_{L^2(0,T;H)}\leq C \quad\text{for all } \delta\in(0,1)\,.
\eeq
Moreover, \eqref{est1}, \eqref{est8_provv} and well-known elliptic regularity results imply that 
\beq
  \label{pier7}
  \delta\l| y_\delta\r|_{L^2(0,T;H^2(\Omega))}\leq C \quad\text{for all } \delta\in(0,1)\,.
\eeq
}%

\subsection{The passage to the limit}

We are now ready to pass to the limit and \pier{conclude the proof of Theorem}~\ref{thm2}. Firstofall, we notice that the estimates
\eqref{est4}, \eqref{est6}, \eqref{est7},  \pier{\eqref{est8_provv} ensure that there exist
\begin{gather}
  \label{pier8}
  y\in H^1(0,T;H)\,,\quad 
  w\in L^2(0,T;V)\,,\quad
  \xi\in L^2(Q)
\end{gather}
such that, at least for a subsequence, 
\begin{gather}
  \label{conv2}
  y_{\delta}\rarrw y \quad\text{in } H^1(0,T;H)\,,\\
    \label{conv2bis}
  \delta \Delta y_{\delta}\rarrw 0 \quad\text{in } L^2(0,T;H),\\
  \label{conv3}
  w_{\delta}\rarrw w \quad\text{in } L^2(0,T;V),\\
  \label{conv4}
  \xi_{\delta}\rarrw\xi \quad\text{in } L^2(Q)
\end{gather}
as $\delta \searrow 0$. Moreover, in view of \eqref{pi} and \eqref{est4} there is some 
$\eta \in H^1(0,T;H) $ such that 
\beq
  \label{conv5}
  \pi (y_{\delta})\rarrw\eta \quad\text{in } \luca{H^1(0,T; H)} \quad \hbox{as } \, \delta \searrow 0.
\eeq
We also note that  $\delta y_{\delta}$ tends to $0$ 
strongly in $L^\infty (0,T; V)$ and weakly in $L^2(0,T;H^2(\Omega))$,
due to \eqref{est4} and \eqref{pier7}.}

\pier{Then, we can pass to the limit in \eqref{1} to readily recover \eqref{1_lim} and in \eqref{pier6}, obtaining~(cf.~\eqref{lambda_delta} too)
 \beq
\label{pier9}
 w = \tau\partial_t y + \xi + \lambda \eta - g  \quad\text{a.e. in } Q.
\eeq
From \eqref{conv2}, using the Ascoli theorem,  it follows that 
\beq
\label{pier9bis}
y_{\delta}\to y \quad\text{in } C^0([0,T];V^*). 
\eeq
\luca{We have to show that}
\beq
\label{pier10}
\luca{\eta = \pi (y)\quad \hbox{and} \quad \xi \in \beta (y) }\quad\text{a.e. in } Q.
\eeq
To this aim, we multiply \eqref{pier6} by the test function $e^{-Rt} y_{\delta}(t) $ and integrate over space and time, with $R>0$ to be chosen soon. As \luca{$ e^{-Rt/2} 
\partial_t y_{\delta} (t) =  \partial_t ( e^{-Rt/2} y_{\delta} (t)) + \frac{R}{2}  e^{-Rt/2} y_{\delta} (t)$}, 
it is easy verify that 
\begin{align}
\frac{\tau}2 \int_\Omega  e^{-RT} | y_{\delta} (T) |^2  
+ \luca{\frac{R}{2}} \tau \int_0^T\!\!\int_\Omega e^{-Rt} \vert y_{\delta} (t) \vert^2 dt 
+  \int_0^T\!\!\int_\Omega  e^{-R t} \, \delta |\nabla y_{\delta}(t)|^2 dt
\nonumber\\
+  \int_0^T\!\!\int_\Omega  e^{-Rt} (\xi_\delta  + \lambda \pi(y_\delta)) (t)\, y_{\delta} (t)\, dt
+  \int_0^T\!\!\int_\Omega  e^{-Rt} g(t) \, y_{\delta} (t)\,  dt
\nonumber\\
= \frac{\tau}2 \Vert y_{0,\delta} \Vert_H^2  
+  \int_0^T  e^{-Rt} \langle y_{\delta}(t) , w_\delta (t) \rangle dt . \label{pier11}
\end{align}
Hence, from \eqref{pier11} we infer that
\begin{align}
&\int_0^T\!\!\int_\Omega  e^{-Rt} \left(\luca{\frac{R}{2}}  \tau \, y_\delta + \lambda \pi(y_\delta) + \xi_\delta\right) 
(t)\, y_{\delta} (t)\, dt 
\nonumber\\
&\leq  \frac{\tau}2 \Vert y_{0,\delta} \Vert_H^2  
+  \int_0^T  e^{-Rt} \langle y_{\delta}(t) , w_\delta (t) \rangle dt 
\nonumber\\
&\quad {}-  \int_0^T\!\!\int_\Omega  e^{-Rt} g(t) \, y_{\delta} (t)\,  dt
- \frac{\tau}2 \int_\Omega  e^{-RT} | y_{\delta} (T) |^2 . 
\label{pier12}
\end{align}
We now choose $R$ in order that the operator $L : L^2(Q) \to L^2(Q)$ defined by
\beq     
\label{defL}
L(v) (x,t) := \luca{\frac{R}{2}} \tau \, v(x,t) + \lambda(x,t)  \pi(v(x,t)) , 
\quad (x,t)\in Q, \quad v\in L^2(Q)
\eeq
be strongly monotone and Lipschitz continuous. In fact, we can take \luca{$R > 2\Vert \lambda \|_{L^\infty (Q)} C_\pi /\tau ,$} where $C_\pi$ stands for a Lipschitz constant for $\pi$ (see \eqref{pi}). Then, it turns out that (see, e.g., \cite[Lemme~2.4, p.~34]{Brezis}) the sum
of $L$ and of the operator induced by $\beta$ on  $L^2(Q)$ (still denoted by $\beta$) is maximal 
monotone. As 
$$
\int_0^T\!\!\int_\Omega  e^{-Rt} v_1 (t) \, v_2 (t) \, dt , \quad v_1, \, v_2 \in L^2(Q),
$$
is an admissible scalar product in $L^2(Q)$, if we can show that
\begin{align}
\limsup_{\delta \searrow 0} \int_0^T\!\!\int_\Omega  e^{-Rt} \left( \luca{\frac{R}{2}} \,\tau \, y_\delta + \lambda \,\pi(y_\delta) + \xi_\delta\right) 
(t)\, y_{\delta} (t)\, dt 
\nonumber \\
\leq \int_0^T\!\!\int_\Omega  e^{-Rt} \left( \luca{\frac{R}{2}} \,\tau \, y + \lambda \, \eta + \xi\right) 
(t)\, y (t)\, dt ,\label{pier13}
\end{align}
then, on account of \eqref{conv2}, \eqref{conv4}, \eqref{conv5} and using a standard result 
for maximal monotone operators (see, e.g., \cite[Prop.~2.5, p.~27]{Brezis}), we actually 
prove \luca{that 
\beq
\label{pier13bis} 
\luca{\frac{R}{2}} \,\tau \, y + \lambda \, \eta + \xi \in (L + \beta) (y) \quad  
\, \text{ in } L^2(Q).
\eeq 
}
In order to check \eqref{pier13}, we pass to the $\limsup $ in the inequality 
\eqref{pier12} noting that 
\begin{align*}
&\lim_{\delta \searrow 0} \frac{\tau}2 \Vert y_{0,\delta} \Vert_H^2  = \frac{\tau}2 \Vert y_{0} \Vert_H^2 \quad \hbox{by \eqref{approx_init}}, \\
&\lim_{\delta \searrow 0}
\int_0^T  e^{-Rt} \langle y_{\delta}(t) , w_\delta (t) \rangle dt = \int_0^T  e^{-Rt} \langle y(t) , w (t) \rangle dt  \quad \hbox{by \eqref{pier9bis} and \eqref{conv3}}, \\
&\lim_{\delta \searrow 0} -  \int_0^T\!\!\int_\Omega  e^{-Rt} g(t) \, y_{\delta} (t)\,  dt
= -  \int_0^T\!\!\int_\Omega  e^{-Rt} g(t) \, y (t)\,  dt \quad \hbox{by \eqref{data_delta} and \eqref{conv2}}, \\
&\limsup_{\delta \searrow 0}
- \frac{\tau}2 \int_\Omega  e^{-RT} | y_{\delta} (T) |^2 
\leq - \liminf_{\delta \searrow 0}
\frac{\tau}2 \int_\Omega  e^{-RT} | y_{\delta} (T) |^2 
\leq - \frac{\tau}2 \int_\Omega  e^{-RT} | y (T) |^2 ,
\end{align*}
where the last inequality is a consequence of the weak convergence of $ y_{\delta} (t)$ to  
$ y (t)$ in $H$ for all $t\in [0,T]$ and of the lower semicontinuity property of the norm in 
$H$. Note that the mentioned weak convergence can also be invoked to obtain the initial condition \eqref{init} from \eqref{init_delta}. Then, using the equality \eqref{pier9} tested by 
$e^{-Rt} y(t) $, we easily conclude that 
\begin{align}
&\limsup_{\delta \searrow 0}
\int_0^T\!\!\int_\Omega  e^{-Rt} \left( \luca{\frac{R}{2}} \tau \, y_\delta + \lambda \pi(y_\delta) + \xi_\delta\right) 
(t)\, y_{\delta} (t)\, dt 
\nonumber\\
&\leq  \frac{\tau}2 \Vert y_{0} \Vert_H^2  
+  \int_0^T  e^{-Rt} \langle y(t) , w (t) \rangle dt 
-  \int_0^T\!\!\int_\Omega  e^{-Rt} g(t) \, y (t)\,  dt
\nonumber\\
&\quad {}
- \frac{\tau}2 \int_\Omega  e^{-RT} | y (T) |^2 =  \int_0^T\!\!\int_\Omega  e^{-Rt} \left( \luca{\frac{R}{2}} \,\tau \, y + \lambda \, \eta + \xi\right) 
(t)\, y (t)\, dt, \nonumber
\end{align}
whence \eqref{pier13} follows.}

\pier{At this point, \luca{we derive the strong convergence}
\beq
\label{pier15}
y_{\delta}(t)\to y (t)\quad\text{in } H \quad \hbox{for all } t\in [0,T]. 
\eeq
To this aim, we take the difference between 
\eqref{pier6} and \eqref{pier9}, multiply
by the test function $e^{-2Rs} (y_{\delta} - y) (s) $ and integrate over $\Omega \times (0,t),$
with $t \in (0,T]. $ Taking into account the already performed computations,
it is straightforward to verify that
\begin{align}
&\frac{\tau}2 \int_\Omega  e^{-2Rt} | (y_{\delta} - y) (t) |^2  
+ \lucabis{\frac R 2\tau \int_0^t e^{-2Rs} \Vert (y_{\delta} - y) (s) \Vert_H^2 ds}
+  \int_0^t\!\!\int_\Omega  e^{-2Rs} \, \delta \, |\nabla y_{\delta}(s)|^2 ds
\nonumber\\
&{}+ \int_0^t\!\!\int_\Omega  e^{-2Rs} \left( \luca{\frac R 2} \,\tau \, y_\delta + \lambda \,\pi(y_\delta) + \xi_\delta - \luca{\frac R 2} \,\tau \, y - \lambda \, \eta - \xi \right) (s)
\,  (y_{\delta} - y) (s)\, ds 
\nonumber\\
&{}= \frac{\tau}2 \int_\Omega  | y_{0,\delta} - y_0  |^2  - \int_0^t\!\!\int_\Omega  e^{-2Rs} \, \delta \Delta y_{\delta}(s) \, y(s) \, ds 
\nonumber\\
&{}\lucabis{+\int_0^t  e^{-2Rs}\left<(y_\delta-y)(s),(w_\delta-w)(s)\right>ds
\quad \hbox{ for all } t\in [0,T].}
\label{pier16}
\end{align}%
Note that the integral on the second line of \eqref{pier16} is non-negative due to the monotonicity of the operator (cf.~\eqref{defL} and \eqref{beta}) $L + \beta$  and to \eqref{incl_delta} and \luca{\eqref{pier13bis}.} 
Moreover, on account of \eqref{approx_init}, \eqref{conv2bis}, \lucabis{\eqref{conv3} and \eqref{pier9bis}} the right hand side of \eqref{pier16} tends to $0$ for all $t\in [0,T]$. 
This implies \eqref{pier15} as well as
\beq
  y_{\delta}\rarr y \quad\text{in } L^2(0,T;H)\,,\quad
  \delta y_{\delta}\rarr 0 \quad\text{in } L^2(0,T;V).\\
\label{pier17}
\eeq
Then, \eqref{conv_y_2} follows from \eqref{pier9bis} and \eqref{pier17}\luca{; moreover, 
the Lipschitz continuity of $\pi$ enables us to conclude that $\pi ( y_{\delta}) 
\rarr \pi (y )$ strongly in $L^2(0,T;H)$, whence $\eta = \pi (y )$
from \eqref{conv5}. At this point, it is not difficult to deduce \eqref{pier10}, and in particular \eqref{incl}, from \eqref{pier13bis}, as
$$\frac{R}{2} \,\tau \, y + \lambda \, \eta \, \hbox{ is the unique element of }\, Ly \, \hbox{ in } \, L^2(Q).
$$ 
Hence, Theorem~\ref{thm2} is completely proved.}
As a further remark, we note that the boundedness of $\{y_\delta \}$ in $L^\infty (0,T; H)$ and the convergence  \eqref{pier15} (or the first one in \eqref{pier17}) ensure that $  y_{\delta}\rarr y $ strongly in $L^p (0,T;H)$ for all $p \in [1, \infty). $}


\section{Continuous dependence for the limit problem}
\setcounter{equation}{0}
\label{cont_dep_lim}

This section is devoted to proving the continuous dependence result for the limit problem. Assume thus 
the hypotheses of Theorem \ref{thm3} and let $(y_i, w_i, \xi_i)$, \pier{$i=1,2$,} be
\luca{any} respective solutions corresponding to the data \pier{in} \eqref{data_dep_lim}.
For simplicity, let us introduce the notation
\begin{gather*}
  y:=y_1-y_2\,, \quad w=w_1-w_2\,, \quad \xi=\xi_1-\xi_2\,,\\
  y_0:=y_0^1-y_0^2\,, \quad g:=g_1-g_2\,, \quad h:=h_1-h_2\,.
\end{gather*}
\pier{Hence,} if we subtract the corresponding equations \eqref{1_lim}--\eqref{2_lim},
for almost every $t\in(0,T)$ we have
\beq
  \label{1_dep_lim}
   \int_{\Omega}{\partial_t y(t)v}+\int_{\Omega}{\nabla w(t)\cdot\nabla v}=\int_\Gamma{h(t)v}
   \quad\text{for all } v\in V\,,
\eeq
\begin{align}
   \label{2_dep_lim}
      \int_\Omega{w(t)v}=\ &\tau\int_{\Omega}{\partial_t y (t)v}
      +\int_\Omega{\xi(t)v} \nonumber\\
      &+\int_\Omega{\lambda(t)\left(\pi(y_1(t))-\pi(y_2(t))\right)v}
      -\int_\Omega{g(t))v}
      \quad\text{for all } v\in V\,.
\end{align}
Please note that hypothesis \eqref{data''_dep_lim} ensure that $(y(t))_\Omega=0$, so that we test \eqref{1_dep_lim} by $v=\n(y(t))$ and \eqref{2_dep_lim} by $v=-y(t)$: summing up, the second and third integral
on the left hand side cancel thanks to \eqref{defN}. \pier{Therefore, for almost every $t\in(0,T)$ we infer that}
\[
\begin{aligned}
  &\left<\partial_ty(t),\n (y(t))\right>+
  \frac{\tau}{2}\frac{d}{dt}\l|y(t)\r|_H^2
  +\int_\Omega{\xi(t)y(t)}\\
  &=\int_\Gamma{h(t)\n(y(t))}
  +\int_\Omega{g(t)y(t)}
  -\int_\Omega{\lambda(t)\left(\pi(y_1(t))-\pi(y_2(t))\right)y(t)}\,.
\end{aligned}
\]
\pier{Then, in view of \eqref{prop1}--\eqref{prop2} and the monotonicity of $\beta$, integrating on $(0,t)$
leads to}
\[
  \begin{aligned}
  \frac{1}{2}\l|y(t)\r|_*^2+\frac{\tau}{2}\l|y(t)\r|_H^2\leq\ &\frac{1}{2}\l|y_0\r|_*^2+\frac{\tau}{2}\l|y_0\r|_H^2
  +\frac{C^2}{2}\l|h\r|_{L^2(\Sigma)}^2+\frac{1}{2}\l|g\r|^2_{L^2(Q)}\\
  &+\frac{1}{2}\int_0^t{\l|y(s)\r|_*^2\,ds}
  +\left(\frac{1}{2}+\l|\lambda\r|_{L^\infty(Q)}C_\pi \right)\int_0^t{\l|y(s)\r|_H^2\,ds}
  \end{aligned}
\]
\pier{for a.e. $t\in (0,T)$ and a certain constant $C>0$. Now, we can apply the Gronwall lemma that implies (updating the value of  $C$)}
\[
  \l|y\r|_{L^\infty(0,T;H)}\leq
  C\left[\l|y_0\r|_H+\l|h\r|_{L^2(\Sigma)}+\l|g\r|_{L^2(Q)}\right]\,,
\]
and the estimate \eqref{dep_lim} is proved.


\section{\pier{Asymptotic} error estimate}
\setcounter{equation}{0}
\label{error_estim}

In this section, we prove the \pier{asymptotic} estimate stated in Theorem \ref{thm4}. To this aim,
 in a first step we deduce an additional uniform estimate to improve the boundedness properties of the solution to the problem \eqref{1}--\eqref{init_delta}.

\subsection{The estimate on $\Delta y_\delta$}

\pier{Here, we want to prove something better than \eqref{est8_provv}. We test \eqref{pier6} by $-\delta^{1/2}\Delta y_\delta(t)$ and deduce that}
  \begin{align}
  &\delta^{3/2}\int_\Omega{\left|\Delta y_\delta(t)\right|^2}-\tau\delta^{1/2}\int_\Omega{\partial_ty_\delta(t)\Delta y_\delta(t)}=
  -\delta^{1/2}\int_\Omega{w_\delta(t)\Delta y_\delta(t)} \nonumber\\
  &+\delta^{1/2}\int_\Omega{\xi_\delta(t)\Delta y_\delta(t)}+
  \delta^{1/2}\int_\Omega{\lambda(t)\pi(y_\delta(t))\Delta y_\delta(t)}
  -\delta^{1/2}\int_\Omega{g(t)\Delta y_\delta(t)}   \label{lap_bis}
  \end{align}
\pier{for almost all $t\in(0,T)$.}
Please note that the second term on the left hand side can be written by integration by parts as
\[
  \frac{\tau}{2}\delta^{1/2}\frac{d}{dt}\int_\Omega{\left|\nabla y_\delta(t)\right|^2}\,;
\]
let us handle the terms on the right hand side. \pier{We integrate by parts, taking into account that $y_\delta$ satisfies Neumann homogeneous boundary conditions. By the Young inequality we have}
\[
  -\delta^{1/2}\int_\Omega{w_\delta(t)\Delta y_\delta(t)}=
  \delta^{1/2}\int_\Omega{\nabla w_\delta(t)\cdot\nabla y_\delta(t)}\leq
  \frac{1}{2}\l|\nabla w_\delta(t)\r|_H^2+\frac{1}{2}\l|\delta^{1/2}\nabla y_\delta(t)\r|_H^2\,,
\]
\[
  -\delta^{1/2}\int_\Omega{g(t)\Delta y_\delta(t)}=
  \delta^{1/2}\int_\Omega{\nabla g(t)\cdot \nabla y_\delta(t)}\leq
  \frac{1}{2}\l|\nabla g(t)\r|_H^2+\frac{1}{2}\l|\delta^{1/2}\nabla y_\delta(t)\r|_H^2\,;
\]
moreover, proceeding formally as we already did in \pier{Subsection~\ref{xi_est}, it is not restrictive to argue with} $\beta_{\pier{\eps}}(y_\delta)$ instead of $\xi_\delta$, so that using monotonicity we \pier{deduce that}
\[
  \delta^{1/2}\int_\Omega{\beta_{\pier{\eps}}(y_\delta(t))\Delta y_\delta(t)}=
  -\delta^{1/2}\int_\Omega{\beta'_{\pier{\eps}}(y_\delta(t))\left|\nabla y_\delta(t)\right|^2}\leq0\,.
\]
Finally, using the hypotheses on $\lambda$ and $\pi$ and the Young inequality, a direct computation leads to
\[
  \begin{split}
  &\delta^{1/2}\int_\Omega{\lambda(t)\pi(y_\delta(t))\Delta y_\delta(t)} \\
  &= {}-\delta^{1/2}\int_\Omega{\pi(y_\delta(t))\nabla\lambda(t)\cdot\nabla y_\delta(t)}
  -\delta^{1/2}\int_{\Omega}{\lambda(t)\pi'(y_\delta(t))\left|\nabla y_\delta(t)\right|^2}\\
  & \leq{}\pier{\frac{1}{2} C_\pi^{\,2}  \l|y_\delta\r|^2_{L^\infty(0,T;H)}\l|\lambda(t)\r|_{W^{1,\infty}(\Omega)}^2}
  +\frac{1}{2}\l|\delta^{1/2}\nabla y_\delta(t)\r|_H^2
  +\l|\lambda\r|_{L^\infty(Q)}C_\pi\delta^{1/2}\int_\Omega{\left|\nabla y_\delta(t)\right|^2}\,.
  \end{split}
\]
Hence, taking all these remarks into account, integrating \eqref{lap_bis} with respect to time, and  
using \pier{\eqref{est4},} \eqref{est6} and the hypotheses \eqref{lambda_delta}--\eqref{data_delta} we obtain
\[
\begin{split}
 & \delta^{3/2}\int_0^t\!\!\int_\Omega{\left|\Delta y_\delta(s)\right|^2\,ds}+\frac{\tau}{2}
 \,\delta^{1/2}\l|\nabla y_\delta(t)\r|_H^2\\
 & \leq
  \frac{\tau}{2}\, \delta^{1/2}\l|\nabla y_{0,\delta}\r|_H^2+C+C\delta^{1/2}\int_0^t{\l|\nabla y_\delta(s)\r|^2_H\,ds}
\end{split}
\]
for \pier{some constant $C>0$; then, using the Gronwall lemma and condition \eqref{approx_init'}, we conclude that
\begin{gather}
  \label{est8}
  \delta^{3/4}\l|\Delta y_\delta\r|_{L^2(0,T;H )} + \delta^{1/4} \l|  y_\delta\r|_{L^2(0,T;V )}  
  \leq C \quad\text{for all } \delta\in(0,1)\,.
\end{gather}
}%

\subsection{Error estimate}
I\pier{n order to prove Theorem~\ref{thm4}, we subtract} \eqref{1} and \eqref{2} to \eqref{1_lim} and \eqref{2_lim}, respectively, and for almost every $t\in(0,T)$ we have
\begin{gather}
  \left<\partial_t(y(t)-y_\delta(t)),v\right>+\int_\Omega{\nabla(w(t)-w_\delta(t))\cdot\nabla v}=0 \quad\text{for all } v\in\H1\,,\\
  \begin{split}
  \int_\Omega{(w(t)-w_\delta(t))v}&=\tau\int_\Omega{\partial_t(y(t)-y_\delta(t))v}-\delta\int_\Omega{\nabla y_\delta(t)\cdot\nabla v}
  +\int_\Omega{(\xi(t)-\xi_\delta(t))v}\\
  &+\int_\Omega{\lambda(t)\left(\pi(y(t))-\pi(y_\delta(t))\right)v}
  \quad\text{for all } v\in\H1\,.
  \end{split}
\end{gather}
Testing the first equation by $\n(y(t)-y_\delta(t))$ and the second by $-(y(t)-y_\delta(t))$\pier{, we sum up 
with the help of \eqref{defN} to cancel the two integrals on the left hand side, 
as usual. Hence,} integrating by parts we obtain for almost every $t\in(0,T)$
\[
  \begin{split}
  \frac{1}{2}\l|y(t)-y_\delta(t)\r|_*^2&+\frac{\tau}{2}\l|y(t)-y_\delta(t)\r|_H^2+\int_0^t\!\!\int_\Omega{(\xi(s)-\xi_\delta(s))(y(s)-y_\delta(s))\,ds}\\
  &=\frac{1}{2}\l|y_0-y_{0,\delta}\r|_*^2+\frac{\tau}{2}\l|y_0-y_{0,\delta}\r|_H^2
  -\delta\int_0^t\!\!\int_\Omega{\Delta y_\delta(s)(y(s)-y_\delta(s))\,ds}\\
  &-\int_0^t\!\!\int_\Omega{\lambda(s)\left(\pi(y(s))-\pi(y_\delta(s))(y(s)-y_\delta(s))\right)\,ds}\,;
  \end{split}
\]
now, the monotonicity of $\beta$, the hypotheses \eqref{pi} and \eqref{lambda_delta},
and the Young inequality lead to
\[
  \begin{split}
  \frac{1}{2}\l|y(t)-y_\delta(t)\r|_*^2&+\frac{\tau}{2}\l|y(t)-y_\delta(t)\r|_H^2\leq
  \frac{1}{2}\l|y_0-y_{0,\delta}\r|_*^2+\frac{\tau}{2}\l|y_0-y_{0,\delta}\r|_H^2\\
  &+\frac{\delta^{1/2}}{2}\l|\delta^{3/4}\Delta y_\delta\r|^2_{L^2(Q)}
  +\left(\frac{1}{2}+\l|\lambda\r|_{L^\infty(Q)}C_\pi\right)\int_0^t{\l|y(s)-y_\delta(s)\r|_H^2\,ds}\,,
  \end{split}
\]
so that \pier{the Gronwall lemma and the estimate \eqref{est8} allow us to infer} that
\[
  \l|y-y_\delta\r|_{L^\infty(0,T;H)}\leq
  C\left[\delta^{1/4}+\l|y_0-y_{0,\delta}\r|_H\right]
\]
for a certain \pier{constant} $C>0$. This finishes the proof. \pier{Let us point out that the error estimate 
\eqref{err} is of order $1/4$ in terms of $\delta$ provided $ \l|y_0-y_{0,\delta}\r|_H \leq C \delta^{1/4}$: 
this condition is ensured for the family defined in \eqref{ellip} whenever, for instance, $y_0 \in V$ (take $z =  
y_{0,\delta} - y_0$ in \eqref{ellip'} and perform the estimate).}


%


\end{document}